\documentclass[12pt,leqno]{article}
\usepackage{amsfonts}
\usepackage{amsmath, amsthm, amsfonts, amssymb, color}
\usepackage{mathrsfs}
\usepackage{color}
\theoremstyle{definition}

\newcommand{\scr}[1]{\mathscr #1}
\definecolor{wco}{rgb}{0.5,0.2,0.3}

\numberwithin{equation}{section} \theoremstyle{remark}

\newcommand{\ua}{\uparrow}

\title{
{\bf  Super Poincar\'e inequality for a dynamic model associated with the two-parameter Dirichlet process}}

\author{
{\bf Weiwei Zhang }\footnote{Postal address: School of Mathematical Sciences, Beijing Normal University, Beijing 100875, China. Email address: 201321130173@mail.bnu.edu.cn.},
 Beijing Normal University }

\begin{document}
\def\tttext#1{{\normalfont\ttfamily#1}}
\def\R{\mathbb R}  \def\ff{\frac} \def\ss{\sqrt} \def\B{\mathbf B}
\def\N{\mathbb N} \def\kk{\kappa} \def\m{{\bf m}}
\def\dd{\delta} \def\DD{\Delta} \def\vv{\varepsilon} \def\rr{\rho}
\def\<{\langle} \def\>{\rangle} \def\GG{\Gamma} \def\gg{\gamma}
  \def\nn{\nabla} \def\pp{\partial} \def\EE{\scr E}
\def\d{\text{\rm{d}}} \def\bb{\beta} \def\aa{\alpha} \def\D{\scr D}
  \def\si{\sigma} \def\ess{\text{\rm{ess}}}
\def\beg{\begin} \def\beq{\begin{equation}}  \def\F{\scr F}
\def\Ric{\text{\rm{Ric}}} \def\Hess{\text{\rm{Hess}}}
\def\e{\text{\rm{e}}} \def\ua{\underline a} \def\OO{\Omega}  \def\oo{\omega}
 \def\tt{\tilde} \def\Ric{\text{\rm{Ric}}}
\def\cut{\text{\rm{cut}}} \def\P{\mathbb P}
\def\C{\scr C}     \def\E{\mathbb E}\def\y{{\bf y}}
\def\Z{\mathbb Z} \def\II{\mathbb I}
  \def\Q{\mathbb Q}  \def\LL{\Lambda}\def\L{\scr L}
  \def\B{\scr B}    \def\ll{\lambda} \def\a{{\bf a}} \def\b{{\bf b}}
\def\vp{\varphi}\def\H{\mathbb H}\def\ee{\mathbf e}\def\x{{\bf x}}
\def\gap{{\rm gap}}\def\PP{\scr P}\def\p{{\mathbf p}}\def\NN{\mathbb N}
\def\cA{\scr A} \def\cQ{\scr Q}\def\cK{\scr K}
\def\LS{C_{LS}}
\maketitle
\begin{abstract} In this paper, we establish the  super Poincar\'e inequality for the two-parameter Dirichlet process when the partition number of the state space is finite. Furthermore, if the partition number is infinite, the super Poincar\'e inequality doesn't hold. To overcome the difficulty caused by the degenerency of the diffusion coefficient on the boundary of the domain, localization method and perturbation argument in \cite{Wbook} are effective.

 \end{abstract} \noindent 2010 Mathematics Subject Classification. Primary:60J70; Secondry:60H10.      \vskip 2cm

\section{Introduction}
The two-parameter Dirichlet process is the natural generalization of the single-parameter Dirichlet process, which first appeared in the context of Bayesian statistics. And both two-parameter Dirichlet process and single-parameter Dirichlet process are pure atomic random measure.

For any $0\leq \alpha< 1$ and $\theta> -\alpha,$ let $\{U_{k}\}_{k\geq 1}$ be a sequence of independent random variables such that $U_{k}$ has $Beta(1-\alpha,\theta+k\alpha)$ distribution. Set
$$V_{1}^{\alpha,\theta}=U_{1},~V_{n}^{\alpha,\theta}=(1-U_{1})\cdots(1-U_{n-1})U_{n},~n\geq 2.$$
The distribution of $(V_{1}^{\alpha,\theta},V_{2}^{\alpha,\theta},\cdots)$ is called two-parameter GEM distribution, denoted by $GEM(\alpha,\theta).$ When $\alpha=0,$ $GEM(0,\theta)$ is the well known GEM ditribution.
Let $P(\alpha,\theta) =(\rho_{1},\rho_{2},\cdots)$ denote $(V_{1}^{\alpha,\theta},V_{2}^{\alpha,\theta},\cdots)$ in descending order, the law of $P(\alpha,\theta)$ is called two-parameter Poisson-Dirichlet distribution. When $\alpha=0,~P(0,\theta)$ is the Poisson-Dirichlet distribution, which was introduced by Kingman in \cite{JCF} to describe the distribution of gene frequecies in a large $k$-th most frequency locus. For $S=\mathbb{N},$ and a sequence of independent identically distributed $S$-value random variables $\{\xi_{i}\}_{i\geq 1}$ with common distribution $\nu_{0}$ on $S,$ which are independent of $P(\alpha,\theta)$,
let $$\Theta_{\alpha,\theta,\nu_{0}}=\sum_{i=1}^{\infty}\rho_{i}\delta_{\xi_{i}},$$ the distribution of $\Theta_{\alpha,\theta,\nu_{0}}$ is called two-parameter Dirichlet process and denoted by $\Pi_{\alpha,\theta,\nu_{0}}.$  Both
$GEM(\alpha,\theta)$ and $PD(\alpha,\theta)$ contain the information only on proportions  while $\Pi_{\alpha,\theta,\nu_{0}}$
contains information on both proportions and types or labels.

In the context of population genetics, both the Poisson-Dirichlet distribution and the Dirichlet process appear as approximations to the equilibrium behavior of certain large populations evolving under the influence of mutation and random genetic drift. That is, the unlabeled infinite-many-neutral-alleles model is time-reversible with the Poisson-Dirichlet distribution $P(0,\theta)$; the labeled infinite-many-neutral-alleles model ( Fleming-Viot process with neutral parent independent mutation ) is time-reversible with the Dirichlet process $\Pi_{0,\theta,\nu_{0}}.$ Based on \cite{FS18}, when the range of parameters is $\alpha=\frac{1}{2}, \theta> -\frac{1}{2},$ the Fleming-Viot process has a two-parameter analogue which is time-reversible with the two-parameter Dirichlet process $\Pi_{\alpha,\theta,\nu_{0}}.$

People care about how the diffusion processes convergent to the stationary distribution. Functional inequalities contain the Poincar\'e  inequality, log-Sobolev inequality, super Poincar\'e inequality and so on ( see \cite{Wbook} and also \cite{Bakry}) can report this information. There are some people have done these work. For example, \cite{FSWX11} established the functional inequalities for the two-parameter extension of the unlabeled infinite-many-neutral-alleles diffusion process.
\cite{S} obtained log-Sobolev inequality for the projection measure of the single-parameter Dirichlet process, and then proved that the Poincar\'e inequality for the single-parameter Dirichlet process holds.
\cite{WZ18} established the super Poincar\'e inequality for the projection measure of the single-parameter Dirichlet process. In this paper, we consider the super Poincar\'e inequality for the projection measure of the two-parameter Dirichlet process.

In the following three subsections, firstly, we briefly recall some facts about the super Poincar\'e inequality; and then introduce the dynamic model we will study; at last, we state the main results of the paper.

\subsection{Super Poincar\'e inequality}
In general, let $(\EE,\D(\EE))$ be a conservative symmetric Dirichlet form on $L^2(\mu)$ for some probability space $(E,\F,\mu)$,   let $(L, \D(L))$ be the associated Dirichlet operator, and let $P_t:=\e^{tL}, t\ge 0 $ be the Markov semigroup.

We say that $(\EE,\mu)$ satisfies the super Poincar\'e inequality with rate function $\bb: (0,\infty)\to (0,\infty),$ if
\beq\label{sp}\mu(f^2)\le r\EE(f,f) +\bb(r) \mu(|f|)^2,\ \ r>0, f\in \D(\EE).\end{equation}
This inequality is equivalent to the uniform integrability of $P_t$, i.e. $P_t$ has zero tail norm:
$$\|P_t\|_{tail}:=\lim_{R\to\infty} \sup_{\mu(f^2)\le 1} \mu((P_tf)^2 1_{\{|P_tf|\ge R\}}) =0,\
 \  t>0.$$ When $P_t$ has a heat kernel with respect to $\mu$, it is also equivalent to the absence of the essential spectrum of $L$ (i.e. the spectrum of  $L$ is purely discrete).
 The super Poincar\'e inequality generalizes the classical Sobolev/Nash type inequalities. For instance, when $\gap(L)>0$, \eqref{sp} with $\bb(r)= \e^{c(1+r^{-1})}$ for some $c>0$ is
 equivalent to the log-Sobolev inequality 
 \beq\label{LS0} \mu(f^2\log f^2)\le \ff 2 C \EE(f,f),\quad  f\in \D(\EE),~\mu(f^2)=1 \end{equation}
 holds for some constant $C>0$; while for a constant $p>1$, \eqref{sp} with $\bb(r)=  c(1+r^{-p})$ holds for some $c>0$ if and only if
 the Nash inequality
 \beq\label{NS}  \mu (f^2)\le C  \EE (f,f)^{\ff p{p+1} }\mu  (|f|)^{\ff 2 {p+1}},\ \ f\in \D(\EE ), \mu (f)=0\end{equation}
 holds for some constant $C>0$, they are also equivalent to
 $$\|P_t-\mu\|_{L^1(\mu)\to L^\infty(\mu)} \le \ff{c'}{(t\land 1)^p} \e^{-\gap(L)t},\quad t>0$$ 
 hold for some constant $c'>0,$
 this implies the semigroup is ultrabounded.

\subsection{Two-parameter dynamic model}
We denote by $\B_{b}(\mathbb{N})$ the set of all bounded Borel measurable functions on
$\mathbb{N},$ $C^{\infty}(\mathbb{R}^{d})$ the set of all infinitely differentiable functions on $\mathbb{R}^{d}$ and $\mathcal{P}_{1}(\mathbb{N})$ the space of
all probability measures on the Borel ¦Ò-algebra $\B(\mathbb{N})$ in $\mathbb{N}.$
For $\varphi\in\B_{b}(\mathbb{N})$ and $\mu\in\mathcal{P}_{1}(\mathbb{N}),$
we denote $\langle\varphi, \mu \rangle =\int_{\mathbb{N}}\varphi d\mu.$
Let $\varphi_{i}\in \B_{b}(\mathbb{N}), 1\leq i\leq d, f\in C^{\infty}(\mathbb{R}^{d})$ and
$$\mathcal{F}_{d}:=\{F|\text{~there exists}~d > 0~\text{such that}~F(\mu)=f(\langle\varphi_{1},¦Ì \rangle,...,\langle\varphi_{d},¦Ì \rangle)\}$$
 for $\mu\in \mathcal{P}_{1}(\mathbb{N}),$  $\mathcal{F}=\bigcup_{d\geq 1}\mathcal{F}_{d}.$
For $x\in \mathbb{N}, F\in \mathcal{F}_{d}$ and $\mu\in \mathcal{P}_{1}(\mathbb{N}),$ we define
\beg{align*}\nabla_{x}F(\mu)&:=\frac{dF}{ds}(\mu+s\delta_{x})\mid_{s=0}\\
&= \sum_{i=1}^{d}\partial_{i}f(\langle\varphi_{1},¦Ì \rangle,\cdots,\langle\varphi_{d},¦Ì \rangle)\varphi_{i}(x).\end{align*}
We write $\nabla F(\mu)$ for the function $x\mapsto \nabla_{x}F(\mu).$

We consider the bilinear form all probability measures on the Borel ¦Ò-algebra $\B(\mathbb{N})$ in $\mathbb{N}.$
\begin{equation} \bigg\{\begin{array}{ll}
\EE(F,G)=\frac{1}{2}\int_{\mathcal{P}_{1}(\mathbb{N})}\langle \nabla F(\mu),\nabla G(\mu)\rangle_{\mu}\Pi_{\alpha,\theta,\nu_{0}}(d\mu), \quad F,G\in\mathcal{F},\\
\mathcal{F}=\{F(\mu)=f(\mu(1),\cdots,\mu(d)):f\in C^{\infty}(\mathbb{R}^{d}),\quad d\geq 1\}.\end{array}\end{equation}
According to \cite[Theorem 2.1]{FS18}, the bilinear form is closable on $L^{2}(\mathcal{P}_{1}(\mathbb{N}),\Pi_{\alpha,\theta,\nu_{0}})$ and its closure $(\EE,\D(\EE))$ is a quasi-regular Dirichlet form. The diffusion process associated with $(\EE,\D(\EE))$ is reversible with the stationary distribution  $\Pi_{\alpha,\theta,\nu_{0}}.$
Denote by $(L,\D(L))$ and the generator of $(\EE,\D(\EE))$ on  $L^{2}(\mathcal{P}_{1}(\mathbb{N}),\Pi_{\alpha,\theta,\nu_{0}}).$
By \cite[Theorem 2.2]{FS18}, $\forall F\in \mathcal{F}_{d},$
we have
\beg{align*}LF(\mu)
&=\frac{1}{2}\sum_{i,j=1}^{d}\mu(i)(\delta_{ij}-\mu(j))\partial_{ij}f(\mu(1),\cdots,\mu(d))\\
&+\frac{1}{2}\sum_{i=1}^{d}\bigg(-\frac{1}{2}-\theta \mu(i)+\frac{1}{2}B_{i}(\mu)\bigg)\partial_{i}f(\mu(1),\cdots,\mu(d))
,\end{align*}
where $$B_{i}(\mu)=\lim_{d\rightarrow\infty}\frac{(d+1)\frac{(\nu_{0}(i))^{2}}{\mu(i)}}{\sum_{i=1}^{d}\frac{(\nu_{0}(i))^{2}}{\mu(i)}+\frac{\nu_{0}(d+1)}{1-\sum_{i=1}^{d}\mu(i)}} $$ exists in $L^{2}(\mathcal{P}_{1}(\mathbb{N}),\prod_{\alpha,\theta,\nu_{0}}).$

Below, we consider the projection case.

For any $d\ge 2,$ we define $$\DD^{ (d)}:=  \{ x \in [0,1]^d:\  \sum_{1\le i\le d}x_i\leq 1 \},~~ x_{d+1}=1-\sum_{1\le i\le d}x_i.$$
Denote $p_i=\nu_{0}(i), 1\leq i\leq d$ and $ p_{d+1}=1-\sum_{i=1}^{d}p_{i}.$
Let \beq\label{projection}\mu^{(d)}(dx)=\ff{\GG(\theta+\frac{d+1}{2})\prod_{i=1}^{d+1}p_{i}}{\pi^{\frac{d}{2}}\GG(\theta+\frac{1}{2})} \ff{\prod_{i=1}^{d+1}x_{i}^{-\frac{3}{2}}}{(\sum_{i=1}^{d+1}\frac{p_{i}^{2}}{x_{i}})^{\theta+\frac{d+1}{2}}}
dx=\rho(x)dx \end{equation}
on  the set $\DD^{ (d)}.$
Define the operator
 \beg{align*}L^{(d)}f(x)
=\frac{1}{2}\sum_{i,j=1}^{d}x_{i}(\delta_{ij}-x_{j})(\partial_{ij}f)(x)+\frac{1}{2}\sum_{i=1}^{d}\bigg(-\frac{1}{2}-\theta x_{i}+\frac{(\theta+\frac{d+1}{2})\frac{p_{i}^{2}}{x_{i}}}{\sum_{i=1}^{d+1}\frac{p_{i}^{2}}{x_{i}}}\bigg)\partial_{i}f(x)
,\end{align*}
$f\in C^{2}(\Delta^{(d)}).$

Dirichlet form
 \beg{align*}\EE^{(d)}(f,g):=\frac{1}{2}\mu^{(d)}\bigg(\sum_{i,j=1}^{d}x_{i}(\delta_{ij}-x_{j})\partial_{i}f\partial_{j}g\bigg)\quad f,g\in \D(\EE^{(d)}),\end{align*}
 with domain $\D(\EE^{(d)})$ being the closure of $C^1(\DD^{(d)})$.
We consider the map
 \beg{align*}
\gamma_{d}:&\mathcal{P}_{1}(\mathbb{N})\rightarrow \Delta^{(d)},\\
&\mu\rightarrow \gamma_{d}(\mu)=(\mu(1),\cdots,\mu(d)).\end{align*}
By \cite[Theorem 3.1]{MAC02}, we have
$$\Pi_{\alpha,\theta,\nu_{0}}\circ\gamma_{d}^{-1}=\mu^{(d)}. $$
That is, $\forall~F, G\in \mathcal{F}_{d},$
we have
\beq\label{eqal1} \Pi_{\alpha,\theta,\nu_{0}}(h(F))=\mu^{(d)}(h(f)), \quad h\in C^{\infty}(\mathbb{R}),\end{equation}
\beg{align*}\label{eqal2}&\int_{\mathcal{P}_{1}(\mathbb{N})}\frac{1}{2}\langle \nabla F(\mu),\nabla G(\mu)\rangle_{\mu}\Pi_{\alpha,\theta,\nu_{0}}(d\mu)\\
&=\frac{1}{2}\int_{\Delta^{(d)}}
\sum_{i,j=1}^{d}x_{i}(\delta_{ij}-x_{j})(\partial_{i}f)(x)(\partial_{j}g)(x)
\mu^{(d)}(dx).\end{align*}
This is an analogue of single-parameter Dirichlet process whose projection on finite partition of $S$ is Dirichlet distribution.

In this paper, we follow the line of thinking in \cite{S} and apply the localization thinking which are effective for these kind questions, and the following are the main results.
\subsection{Main results}
Let $\nu_{0}$ be the probability measure on type space $S=\mathbb{N}$ in the definition of the two parameter Dirichlet process. Denote $d:=\sharp\{i\in S, \nu_{0}(i)> 0\}.$
 \beg{thm}\label{d_{1}}
If $d<\infty,$ then there exists a constant $c> 0$ such that
 the super Poincar\'e inequality
 $$\Pi_{\alpha,\theta,\nu_{0}}(F^2)\le 1+r\EE(F,F) + c(1+r^{-\frac{1}{2}((\theta+\frac{d}{2})(2d+1)-1)})\Pi_{\alpha,\theta,\nu_{0}}(|F|)^2,\quad r>0, F\in \mathcal{F}_{d}$$ holds,
 where $\mathcal{F}_{d}$ is defined in \eqref{finite}.
 \end{thm}

 \beg{thm}\label{d_{2}}
If $d=\infty,$ then
 the super Poincar\'e inequality doesn't hold.
  \end{thm}
We remark that there is a question which we haven't finished: if $d=\infty,$ does
the Poincar\'e inequality for $\Pi_{\alpha,\theta,\nu_{0}}$ hold?

To establish the super  Poincar\'e inequality for the measured-value process, we firstly establish the super  Poincar\'e inequality for the projection measure of $\Pi_{\alpha,\theta,\nu_{0}}$ in Section 2, then we prove Theorem \ref{d_{1}} and Theorem \ref{d_{2}} in Section 3.

\section{Super Poincar\'e inequality for $\mu^{(d)}$}
To establish the super Poincar\'e inequality for the projection measure of $\Pi_{\alpha,\theta,\nu_{0}}$ with an explicit rate function $\bb$, the main difficulty   comes from the degeneracy of the diffusion coefficient on the boundary $$\pp \DD^{(d)}=\Big\{x=(x_i)_{1\le i\le d}\in \DD^{(d)}:\ \min\{x_i: 1\le i\le d+1\}=0\Big\},\ \ x_{d+1}:=1-\sum_{i=1}^d x_i.$$

We have two methods to establish the super Poincar\'e inequality. firstly, from  \cite{WZ18}, we have known the super Poincar\'e inequality for another probability measure $\widetilde{\mu}^{(d)}$~already, then we can get the local super Poincar\'e inequality for the measure $\mu^{(d)},$
 by \cite[Theorem 2.1,Lemma 3.4]{WZ18}, we can establish the super Poincar\'e inequality for $\mu^{(d)}.$ Secondly, from \cite{WZ18}, we have known the super Poincar\'e inequality for another probability measure $\widetilde{\mu}^{(d)}$~already,  by the perturbation result for the super Poincar\'e inequality \cite[Theorem 3.4.7]{Wbook}, we can establish the super Poincar\'e inequality for $\mu^{(d)}.$

\subsection{Preparations}
$ Assumption (A)$:
Let $(E,\F,\mu)$ be a separable complete probability space, and let
$(\EE,\D(\EE))$ be a conservative symmetric local Dirichlet form on
$L^2(\mu)$ as the closure of
$$\EE(f,g)=\mu(\GG(f,g)),\ \ f,g\in \D_0(\GG),$$
where $\GG: \D(\GG)\times \D(\GG)\to \scr B(E)$ is a
  positive definite symmetric bilinear mapping,   $\scr B(E)$ is the set of all
$\mu$-a.e. finite measurable real functions on $E$,   $\D(\GG)$ is a sub-algebra of   $\B(E)$, and $ \D_0(\GG):=\{f\in \D(\GG): f^2,\GG(f,f)\in L^1(\mu)\}$ such that
\beg{enumerate} \item[(a)]  $\D_0(\GG)$ is dense in $L^2(\mu)$.
\item[(b)] $\D(\GG)$ is closed under combinations with   $\psi\in C([-\infty,\infty])$ such that $\psi$ is $C^1$ in $\R$ and  $\psi'$ has compact support, and $\GG(\psi\circ f, g)= \psi'(f) \GG(f,g)\ \mu$-a.e.  for $f,g\in \D(\GG)$.
\item[(c)] $\GG(fg,h)= g\GG(f,h)+ f\GG(g,h)\ \mu$-a.e. for $f,g,h\in \D(\GG)$.  \end{enumerate}
  Let $\phi\in \D(\GG)$ be an unbounded  nonnegative   function and let
  $$h(s):=\sup_{\phi\leq s}\Gamma(\phi,\phi), s\geq s_{0}>0. $$
   \beq\label{HU}  D_s:=\{\phi\le s\},\ \ s\ge s_{0}>0,\end{equation}   where $\sup_{\emptyset}=0$ by convention.
    \beq\label{D2} 0<\ll(s):= \inf\{\EE(f,f): \mu(f^2)=1, f|_{D_s}=0\}\uparrow \infty\ \text{as}\ s\uparrow\infty.\end{equation}
To verify Assumption $(A),$ we can take $$\D(\GG)=\Big\{f\in C(\DD^{(d)}; [-\infty,\infty]): f \ \text{is\ finite\ and}\  C^1 \ \text{in}\ \DD^{(d)}\setminus \{x_{i}=0,~1\leq i\leq d+1\}\Big\},$$ and let
$$\GG(f,g)(x)= 1_{\{x_{k}>0,~1\leq k\leq d\}}\sum_{i,j=1}^d x_i(\delta_{ij}-x_{j})(\pp_i f)(x) (\pp_j g)(x),\ \ f,g\in \D(\GG).$$
Obviously, conditions (a)-(c) in $ Assumption (A)$ hold.

We set
\beq\label{PHI}\phi(x)= \sum_{1\leq i\leq d+1}x_{i}^{-2},\ \ x=(x_i)_{1\le i\le d}\in \DD^{(d)},\end{equation}
then for $s\geq s_{0}>0,$
\beq\label{DS}D_s:=\{\phi\le s\}= \big\{x\in \DD^{(d)}:\sum_{1\leq i\leq d+1}x_{i}^{-2}\leq s\big\},\end{equation}
\beg{align*}h(s)
&:=\sup_{\phi\leq s}\Gamma(\phi,\phi)\\
&=\sup_{\phi\leq s}\sum_{i,j=1}^d x_i(\delta_{ij}-x_{j})(\pp_i\phi)(x) (\pp_j\phi)(x)\\
&=\sup_{\phi\leq s}\sum_{i,j=1}^d x_i(\delta_{ij}-x_{j})\bigg(\frac{-2}{x_{i}^{3}}+\frac{2}{x_{d+1}^{3}}\bigg) \bigg(\frac{-2}{x_{j}^{3}}+\frac{2}{x_{d+1}^{3}}\bigg)\\
&\leq c_{3}s^{\frac{5}{2}}.\end{align*}
In the following subsection, we estimate $\ll(s).$

\subsection{Estimate on $\ll(s)$}
Let $\ll(s)= \inf\{\EE^{(d)}(f,f):\ f\in C^1(\DD^{(d)}), \mu^{(d)}(f^2)=1, f|_{D_s}=0\}$.  We will adopt the following Cheeger inequality to estimate $\ll(s)$. Let $$\pp D_s= \{x\in \DD^{(d)}: \exists 1\leq i\leq d+1,x_{i}=b_{i}s^{-\frac{1}{2}}\},\ s\ge s_{0}> 0, b_{i}> 0.$$

 \beg{lem}\label{L2.2}  If there exists a function $f\in C^2(\DD^{(N)}\setminus \{x_{i}=0,~1\leq i\leq d+1\})$ such that
 \beq\label{GG} |\sigma\nabla f(x)|:= \sqrt{\sum_{i=1}^{d}\bigg(\sum_{j=1}^d x_i(\delta_{ij}-x_{j})(\pp_j f)(x)\bigg)^{2}}\le a_1,\ \ |L^{(d)}f(x)|\ge a_2, \ \ x\in D_s^c\end{equation}
 holds for some constants $a_1,a_2>0$, and that
 \beq\label{GG2} \lim_{r\to\infty} \int_{x\in\pp D_r} \bigg[\sum_{i=1}^d \bigg(\sum_{j=1}^d x_{i} (\delta_{ij}-x_j)\pp_j f(x)\bigg)^{2}\bigg]^{\frac{1}{2}}\rho(x)dA=0.\end{equation}
 Then
 $$\ll(s)\ge \ff{a_2^2}{4a_1}.$$\end{lem}

 \beg{proof} By \eqref{GG}, we assume that $L^{(d)}f|_{D_s^c}\ge a_2$, otherwise simply use $-f$  replacing $f$.
 Let $\si(x)=\{x_i(\delta_{ij}-x_{j})\}_{1\le i,j\le d}$.
For any nonnegative $g\in C^1(\DD^{(d)})$ with $g|_{D_s}=0$, we have $g|_{\pp D_s}=0$, so that by integration by parts formula,
\beq\label{GG3} \beg{split} &a_2 \mu^{(d)} (g) \le \mu^{(d)} (g L^{(d)}f)\\
& =\lim_{r\to \infty} \int_{D_r\setminus D_s} (\rr g L^{(d)}f)(x)\d x\\
&\le -\mu^{(d)}\big(\<\si \nn g, \si\nn f\>\big) \\
&+\|g\|_\infty\limsup_{r\to\infty} \int_{\pp D_r} \bigg[\sum_{i=1}^d \bigg(\sum_{j=1}^dx_i(\delta_{ij}-x_{j})\pp_j f(x)\bigg)^{2}\bigg]^{\frac{1}{2}} \rr(x)  \d A,\end{split}\end{equation}
where $A$ is the area measure on $\pp D_r$ induced by the Lebesgue measure.
Combining this with \eqref{GG}, \eqref{GG2} and \eqref{GG3}, we obtain
$$a_2 \mu^{(d)} (g)\le |\mu^{(d)}\big(\<\si \nn g, \si\nn f\>\big)|\le \ss{a_1 } \mu^{(d)}(|\si\nn g|).$$ Therefore, for any $g\in C^1(\DD^{(d)})$ with $g|_{D_s}=0$,
$$\mu^{(d)}(g^2)\le \ff{\ss{a_1}}{a_2} \mu^{(d)}(|\si\nn g^2|)\le \ff{2\ss{a_1}}{a_2}\ss{\mu^{(d)}(g^2)\mu^{(d)}(|\si\nn g|^2)}.$$
Noting that $\mu^{(d)}(|\si\nn g|^2)=\EE^{(d)}(g,g)$, we arrive at
$$\mu^{(d)}(g^2)\le \ff{4a_1}{a_2^2}\EE^{(d)}(g,g),\ \ g\in C_b^1(\DD^{(d)}), g|_{D_s}=0. $$
This finishes the proof. \end{proof}

\beg{lem}\label{L2.3}For the operator $L^{(d)},$ there exist  constants $s_0, c_6>0$ such that
$$\ll(s)\ge c_{6}s^{\frac{7}{8}},\ \ s\ge s_0> 0.$$\end{lem}
\beg{proof}  We take
$$f(x)=\sum_{1\leq i\leq d} x_{i}^{\frac{1}{4}}, \ \ x\in \DD^{(d)}.$$
 Then $\forall x\in D_{s}^{c},$ we have
 \beq\beg{split} &|\sigma\nabla f(x)|:= \sqrt{\sum_{i=1}^{d}\bigg(\sum_{j=1}^d x_i(\delta_{ij}-x_{j})(\pp_j f)(x)\bigg)^{2}}\\
 &\le \bigg[\sum_{i=1}^{d}[\frac{1}{4}x_{i}^{\frac{1}{4}}-\frac{1}{4}x_{i}\sum_{j=1}^{d}x_{j}^{\frac{1}{4}}]^{2}\bigg]^{\frac{1}{2}}
\leq c_{4}s^{\frac{-1}{8}},\end{split}\end{equation}
and
\beq\beg{split} &|L^{(d)}f(x)| \\ &=\bigg|\sum_{i,j=1}^{d}x_{i}(\delta_{ij}-x_{j})(\partial_{ij}f)(x)+\sum_{i=1}^{d}\bigg(-\frac{1}{2}-\theta x_{i}+\frac{(\theta+\frac{d+1}{2})\frac{p_{i}^{2}}{x_{i}}}{\sum_{i=1}^{d+1}\frac{p_{i}^{2}}{x_{i}}}\bigg)\partial_{i}f(x)\bigg|\\
&\geq c_{5}s^{\frac{3}{8}}.\end{split}\end{equation}
Denote $$\pp D_{r,j}=\{x|x_{j}=b_{j}r^{-\frac{1}{2}}\},$$
we get
\beq\label{bianjie} \beg{split}
&\int_{\pp D_r} \ff{\prod_{i=1}^{d+1}x_{i}^{-\frac{3}{2}}}{(\sum_{i=1}^{d+1}\frac{p_{i}^{2}}{x_{i}})^{\theta+\frac{d+1}{2}}}\d A \\
&=\sum_{j=1}^{d+1}\int_{\pp D_{r,j}} \ff{\prod_{i\neq j, 1\leq i\leq d+1}x_{i}^{-\frac{3}{2}}}{(\sum_{i\neq j, 1\leq i\leq d+1}\frac{p_{i}^{2}}{x_{i}})^{\theta+\frac{d+1}{2}}}\d A\\
& =\sum_{i=1}^{d}(1-\sum_{j\neq i}b_{j}r^{-\frac{1}{2}})^{\theta+\frac{1}{2}-d} \int_{\DD^{(d-1)}}\frac{\Big(1-\sum_{1\le i\le d-1}x_i\Big)^{-\frac{3}{2}} \prod_{i=1}^{d-1} x_i^{-\frac{3}{2}} }{(\sum_{1\leq i\leq d}\frac{p_{i}^{2}}{x_{i}})^{\theta+\frac{d+1}{2}}}\d x\end{split}\end{equation}
is bounded,
so $$\limsup_{r\to\infty} \int_{\pp D_r} \sum_{i=1}^d \bigg(\sum_{j=1}^dx_i(\delta_{ij}-x_{j})\pp_j f(x)\bigg)^{2} \rr(x)  \d A=0.$$
We derive from   Lemma \ref{L2.2} that
$$\ll(s)\ge \ff{c_{5}^{2}s^{\frac{3}{4}}}{4c_{4}s^{\frac{-1}{8}}}=c_{6}s^{\frac{7}{8}} ,\ \ s\ge s_0> 0.$$
 \end{proof}

\subsection{Localization method}
\beg{thm}\label{super Poincare 1}
Let $\mu^{(d)}$ defined as \eqref{projection}, then the super Poincar\'e inequality
 \beq\mu^{(d)} (f^2)\le r\EE^{(d)}(f,f)+\beta(r)\mu^{(d)}(|f|)^2,\ \ r>0, f\in \D(\EE^{(d)})\end{equation}
 holds with $\beta(r)= c_{13}(1+r^{-\{[(2\theta+d)d+(\theta+\frac{d}{2}-1)]+\frac{3}{7}\}}),~c_{13}> 0.$
\end{thm}

\begin{proof}
We know
$$\mu^{(d)}=\ff{\GG(\theta+\frac{d+1}{2})\prod_{i=1}^{d+1}p_{i}}{\pi^{\frac{d}{2}}\GG(\theta+\frac{1}{2})} \ff{\prod_{i=1}^{d+1}x_{i}^{-\frac{3}{2}}}{(\sum_{i=1}^{d+1}\frac{p_{i}^{2}}{x_{i}})^{\theta+\frac{d+1}{2}}}
dx.$$
Denote
$$\widetilde{\mu}^{(d)}=\ff{\GG(|\aa|_1)}{\prod_{1\le i\le d+1} \GG(\aa_i)} (1-|x|_1)^{\aa_{d+1}-1}\prod_{1\le i\le d} x_i^{\aa_i-1}dx.$$
We set $$\widetilde{\mu}^{(d)}=\mu^{(d)}(e^{W}), $$
so $$e^{W}=\ff{\pi^{\frac{d}{2}}\GG(\theta+\frac{1}{2})\GG(|\aa|_1)(\sum_{i=1}^{d+1}\frac{p_{i}^{2}}{x_{i}})^{\theta+\frac{d+1}{2}}}{\GG(\theta+\frac{d+1}{2})\prod_{i=1}^{d+1}p_{i}\prod_{1\le i\le d+1} \GG(\aa_i)}\prod_{1\le i\le d+1} x_i^{\aa_i+\frac{1}{2}}dx.$$
When $$\alpha_{i}\geq\theta+\frac{d}{2}, 1\leq i\leq d+1,$$
 there are constants $C_{1}$ and $C_{2}$ such that
$$C_{1}s^{-\sum_{i=1}^{d+1}(\alpha_{i}-(\theta+\frac{d}{2}))}\leq e^{W}\leq C_{2}.$$
So the local super Poincar\'e inequality becomes
 \beg{align*}
&\mu^{(d)}(f^2)\le C_{3}rs^{-\sum_{i=1}^{d+1}(\alpha_{i}-(\theta+\frac{d}{2}))}\EE^{(d)}(f,f)+ C_{3}\bb_s(r) s^{-\sum_{i=1}^{d+1}(\alpha_{i}-(\theta+\frac{d}{2}))}\mu^{(d)}(|f|)^2\\
&=C_{3}v\EE^{(d)}(f,f)+ C_{3}\bb_s(v)\mu^{(d)}(|f|)^2, r>0, f\in \D(\EE), f|_{D_s^c}=0,
\end{align*}
where $\bb_s(v)=c_{14}(1+v^{-p'}).$
When  $\alpha_{i}=\theta+\frac{d}{2}, 1\leq i\leq d+1,$ we get the smallest $p'=[(2\theta+d)d+(\theta+\frac{d}{2}-1)].$
By \cite[Theorem 2.1]{WZ18} without the condition that $h(s)<\infty,$ we know the super Poincar\'e inequality holds with
$$\beta(r)=c_{13}(1+r^{-(p'+\frac{3}{7})}).$$
\end{proof}

\subsection{Perturbation method}

Below is the perturbation theorem which is similar to \cite[Theorem 3.4.7]{Wbook}.

\beg{thm}\label{supoin}Under Assumption $(A),$ let~$W$~is bounded on~$\{\phi\leq r\}$~for any ~$r> 0.$ Let $S(W)\in\B$ be such that for any nonnegative $f\in\D(\Gamma)$ with $suppf\subset\{\phi\leq r\}$ for some $r> 0,$ one has
    $$\int\Gamma(f,W)d\mu\geq-\int fS(W)d\mu\in\mathbb{R}.$$
    Put
    $$\varphi(r)=\sup\{e^{W}:\phi\leq r\}, $$
    $$\psi(r)=\frac{1}{4}\sup\{\Gamma(W,W)+2S(W):\phi\leq r\}.$$
    If there exist $c_{1},p> 0$ such that
   \beq\label{rao}\mu(f^{2}e^{W})\leq r\mu(e^{W}\Gamma(f,f))+c_{1}(1+r^{-p})\mu(e^{W}|f|)^{2},\end{equation}
    then \eqref{sp} holds with
 \beg{align*}
&\beta(r)=\frac{c_{1}}{1-\bigg\{\frac{2(1+r)^2h(2s)}{\lambda(2s)s^2}+\frac{r(1-\frac{2(1+r)^2h(2s)}{\lambda(2s)s^2})[(1+r)^2s^{-2}h(3s)+\psi(3s)]}{2+r[(1+r)^2s^{-2}h(3s)+\psi(3s)]} \bigg\}}\\
&\cdot\bigg(1+\bigg(\frac{2+r[(1+r)^2s^{-2} h(3s)+\psi(3s)]}{r(1-\frac{2(1+r)^2h(2s)}{\lambda(2s)s^2})}\bigg)^{p}\bigg)\varphi(3s) ,\end{align*}
 where $s=c_{2}r^{\frac{-8}{7}}.$ \end{thm}

\beg{proof}
Let $S(W)\in\B$ be such that for any nonnegative $f\in\D(\Gamma)$ with $suppf\subset\{\phi\leq r\}$ for some $r> 0,$ one has
    $$\int\Gamma(f,W)d\mu\geq-\int fS(W)d\mu\in\mathbb{R}.$$
    Put
    $$\varphi(r)=\sup\{e^{W}:\phi\leq r\}, $$
    $$\psi(r)=\frac{1}{4}\sup\{\Gamma(W,W)+2S(W):\phi\leq r\}.$$
It suffices to consider  $f\in \D_0(\GG)$. For any $s\ge s_0$ and small $\vv\in (0,1)$, let
$\varphi_i\in C^1([0,\infty])$ with $0\le \varphi_i\le 1, |\varphi_i'(s)|\le (1+\vv) s^{-1}, i=1,2$ such that
$$\varphi_1|_{[0, s]}=0,\ \ \varphi_1|_{[2s,\infty]}=1;\ \ \varphi_2|_{[0,2s]}=1,\ \ \varphi_2|_{[3s, \infty]}=0.$$
Let $f_i= f \cdot \varphi_i\circ\phi, 1\le i\le 2.$ Then
$f^2\le f_1^2+f_2^2$ and by conditions (b) and (c),
\beg{align*} \GG(f_1,f_1)&\le 2\GG(f,f)+ 2(1+\vv)^2 f^2 s^{-2}h(2s)\\
& \le  2\GG(f,f) + \ff{2(1+\vv)^2h(2s)}{s^2} f^2,\\
 \GG(f_2,f_2)&\le 2 \GG(f,f) + \ff{2(1+\vv)^2h(3s)}{s^2} f^2.\end{align*} In particular, $f_1,f_2\in \D_0(\GG)\subset \D(\EE)$.
By the definition, we obtain
\beq\label{wai}\mu(f_{1}^{2})\leq\frac{1}{\lambda(s)}\mu\bigg(2\Gamma(f,f)+\ff{2(1+\vv)^2h(2s)}{s^2} f^2\bigg) . \end{equation}
We choose $f_{2}\exp\{-\frac{W}{2}\}$ as the test function in \eqref{rao}, so
we get that
\beq\beg{split}\label{limian}
&\mu(f_{2}^{2})
\leq r\mu\bigg(\Gamma(f_{2},f_{2})-\frac{1}{2}\Gamma((f_{2})^{2},W)+\frac{1}{4}(f_{2})^{2}\Gamma(W,W)    \bigg)\\
&+c_{1}(1+r^{-p})\varphi(3s)\mu(|f|)^{2}\\
&\leq 2r\mu (\GG(f,f)+ (1+\vv)^2 f^2 s^{-2}h(3s))+\frac{r}{4}\mu((f_{2})^{2} \{\Gamma(W,W)+2S(W)\})\\
&+c_{1}(1+r^{-p})\varphi(3s)\mu(|f|)^{2}\\
&\leq 2r\mu(\Gamma(f,f))+r[2(1+\vv)^2s^{-2} h(3s)+\psi(3s)]\mu(f^{2})+c_{1}(1+r^{-p})\varphi(3s)\mu(|f|)^{2}.
\end{split}\end{equation}
Put \eqref{limian} and \eqref{wai} together, then

\beg{align*}&\mu(f^{2})\leq \mu(f_{1}^{2})+\mu(f_{2}^{2})\\
&\leq
\frac{1}{\lambda(s)}\mu\bigg(2\Gamma(f,f)+\ff{2(1+\vv)^2h(2s)}{s^2} f^2\bigg)\\
& + 2r\mu(\Gamma(f,f))+r[(1+\vv)^2s^{-2}h(3s)+\psi(3s)]\mu(f^{2})+c_{1}(1+r^{-p})\varphi(3s)\mu(|f|)^{2}\\
&\leq (\frac{2}{\lambda(s)}+2r )\mu(\Gamma(f,f))+\bigg\{\frac{2(1+\vv)^2h(2s)}{\lambda(2s)s^2}+r[2(1+\vv)^2s^{-2} h(3s)+\psi(3s)]\bigg\}\mu(f^{2})\\
&+c_{1}(1+r^{-p})\varphi(3s)\mu(|f|)^{2}.
 \end{align*}
Set
$$\varepsilon=\frac{2r+\frac{2}{\lambda(s)}}{1-\bigg\{\frac{2(1+\vv)^2h(2s)}{\lambda(2s)s^2}+r[2(1+\vv)^2s^{-2} h(3s)+\psi(3s)]\bigg\}},$$
and choose $$\lambda(s)=r^{-1},$$
so
$$r=\frac{\varepsilon(1-\frac{2(1+\vv)^2h(2s)}{\lambda(2s)s^2})}{4+\varepsilon[(1+\vv)^2s^{-2} h(3s)+\psi(3s)]}. $$
Thus, \beg{align*}
&\beta(\varepsilon)=\frac{c_{1}}{1-\bigg\{\frac{2(1+\vv)^2h(2s)}{\lambda(2s)s^2}+\frac{\varepsilon(1-\frac{2(1+\vv)^2h(2s)}{\lambda(2s)s^2})[(1+\vv)^2s^{-2}h(3s)+\psi(3s)]}{4+\varepsilon[(1+\vv)^2s^{-2}h(3s)+\psi(3s)]} \bigg\}}\\
&\cdot\bigg(1+\bigg(\frac{4+\varepsilon[(1+\vv)^2s^{-2} h(3s)+\psi(3s)]}{\varepsilon(1-\frac{2(1+\vv)^2h(2s)}{\lambda(2s)s^2})}\bigg)^{p}\bigg)\varphi(3s) ,\end{align*}
where $s=c_{2}\varepsilon^{\frac{-8}{7}}.$
\end{proof}

\beg{thm}\label{super Poincare}
Let $\mu^{(d)}$ defined as \eqref{projection}, then the super Poincar\'e inequality
 \beq\label{EXTE}\mu^{(d)} (f^2)\le r\EE^{(d)}(f,f)+\beta(r)\mu^{(d)}(|f|)^2,\ \ r>0, f\in \D(\EE^{(d)})\end{equation}
 holds, where $\beta(r)= c_{12}(1+r^{-\frac{1}{2}((\theta+\frac{d}{2})(2d+1)-1)}).$
\end{thm}
\begin{proof}
We know
$$\mu^{(d)}=\ff{\GG(\theta+\frac{d+1}{2})\prod_{i=1}^{d+1}p_{i}}{\pi^{\frac{d}{2}}\GG(\theta+\frac{1}{2})} \ff{\prod_{i=1}^{d+1}x_{i}^{-\frac{3}{2}}}{(\sum_{i=1}^{d+1}\frac{p_{i}^{2}}{x_{i}})^{\theta+\frac{d+1}{2}}}
dx.$$
Denote
$$\widetilde{\mu}^{(d)}=\ff{\GG(|\aa|_1)}{\prod_{1\le i\le d+1} \GG(\aa_i)} (1-|x|_1)^{\aa_{d+1}-1}\prod_{1\le i\le d} x_i^{\aa_i-1}dx.$$
We set $$\widetilde{\mu}^{(d)}=\mu^{(d)}(e^{W}), $$
so $$e^{W}=\ff{\pi^{\frac{d}{2}}\GG(\theta+\frac{1}{2})\GG(|\aa|_1)(\sum_{i=1}^{d+1}\frac{p_{i}^{2}}{x_{i}})^{\theta+\frac{d+1}{2}}}{\GG(\theta+\frac{d+1}{2})\prod_{i=1}^{d+1}p_{i}\prod_{1\le i\le d+1} \GG(\aa_i)}\prod_{1\le i\le d+1} x_i^{\aa_i+\frac{1}{2}}dx.$$
 Let $S(W)\in\B$ be such that for any nonnegative $f\in\D(\Gamma)$ with $suppf\subset\{\phi\leq r\}$ for some $r> 0,$ one has
    $$\int\Gamma(f,W)d\mu\geq-\int fS(W)d\mu\in\mathbb{R}.$$
    Put
    $$\varphi(r)=\sup\{e^{W}:\phi\leq r\}, $$
    $$\psi(r)=\frac{1}{4}\sup\{\Gamma(W,W)+2S(W):\phi\leq r\}.$$
 From \cite[Theorem 1.1]{WZ18}, the super Poincar\'e inequality holds for $\widetilde{\mu}^{(d)},$
 $$\mu^{(d)}(f^{2}e^{W})\leq r\mu^{(d)}(e^{W}\Gamma(f,f))+c_{1}(1+r^{-p})\mu^{(d)}(e^{W}|f|), $$
 where $p=\Sigma_{i=1}^{d}1\vee(2\alpha_{i})+(\alpha_{d+1}-1)^{+}, c_{1}> 0.$
Then from Theorem \ref{supoin}, we know
\beg{align*}
&\beta(\varepsilon)=\frac{c_{1}}{1-\bigg\{\frac{2(1+\vv)^2h(2s)}{\lambda(2s)s^2}+\frac{\varepsilon(1-\frac{2(1+\vv)^2h(2s)}{\lambda(2s)s^2})[(1+\vv)^2s^{-2}h(3s)+\psi(3s)]}{2+\varepsilon[(1+\vv)^2s^{-2} h(3s)+\psi(3s)]} \bigg\}}\\
&\cdot\bigg(1+\bigg(\frac{2+\varepsilon[(1+\vv)^2s^{-2} h(3s)+\psi(3s)]}{\varepsilon(1-\frac{2(1+\vv)^2h(2s)}{\lambda(2s)s^2})}\bigg)^{p}\bigg)\varphi(3s) .\end{align*}
When $s\ge s_0,$ we have
$$h(s):=\sup_{\phi\leq s}\Gamma(\phi,\phi)\leq c_{7}s^{\frac{5}{2}},$$
$$\ll(s)\ge c_{6}s^{\frac{7}{8}} ,$$
$$\frac{2(1+\vv)^2h(2s)}{\lambda(2s)s^2}\leq c_{8}s^{\frac{-3}{8}},$$
$$\psi(3s)\leq c_{9}s^{-1}.$$
Thus, $$\frac{\varepsilon(1-\frac{2(1+\vv)^2h(2s)}{\lambda(2s)s^2})[(1+\vv)^2s^{-2}h(3s)+\psi(3s)]}{2+\varepsilon[(1+\vv)^2s^{-2} h(3s)+\psi(3s)]}\leq c_{10}s^{-1},$$
$$\frac{2+\varepsilon[(1+\vv)^2s^{-2} h(3s)+\psi(3s)]}{\varepsilon(1-\frac{2(1+\vv)^2h(2s)}{\lambda(2s)s^2})}\leq c_{11}\varepsilon^{-1}.$$
When $\alpha_{i}\geq \theta+\frac{d}{2},$
$\varphi(3s)$ is bounded.
So
 $$\bb(r)=c_{12}(1+r^{-\frac{1}{2}\{\Sigma_{i=1}^{d}1\vee(2\alpha_{i})+(\alpha_{d+1}-1)^{+}\}}).$$
When $\alpha_{i}=\theta+\frac{d}{2}, 1\leq i\leq d+1,$ we have $$r^{-\{\sum_{i=1}^{d}1\vee(2\alpha_{i})+(\alpha_{d+1}-1)^{+}\}} =r^{-((2\theta+d)d+\theta+\frac{d}{2}-1)}.$$
So we get
$$\mu^{(d)} (f^2)\le r\EE^{(d)}(f,f) + c_{12}(1+r^{-((\theta+\frac{d}{2})(2d+1)-1)})\mu^{(d)}(|f|)^2,\ f\in C^1(\DD^{(d)}).$$
\end{proof}

\section{Proof of Theorem \ref{d_{1}} and \ref{d_{2}} }
\subsection{Proof of Theorem \ref{d_{1}}}
\beg{proof}
For $d\geq 2, \forall F\in \mathcal{F}_{d},$ which is defined in \eqref{finite}.
As $$(2\theta+d)d+\theta+\frac{d}{2}-1+\frac{3}{7}\geq (2\theta+d)d+\theta+\frac{d}{2}-1, $$
by Theorem \ref{super Poincare 1} and Theorem \ref{super Poincare},
we have the super Poincar\'e inequality for $\mu^{(d)}$ and $\EE^{(d)},$ 
$$\mu^{(d)} (f^2)\le r\EE^{(d)}(f,f) + c(1+r^{-((2\theta+d)d+\theta+\frac{d}{2}-1)})\mu^{(d)}(|f|)^2,\ \ r>0, f\in \D(\EE^{(d)}).$$
From \eqref{eqal1}, we have the super Poincar\'e inequality for $\Pi_{\alpha,\theta,\nu_{0}}$ and $\EE,$ 
$$\Pi_{\alpha,\theta,\nu_{0}}(F^2)\le r\EE(F,F) + c(1+r^{-((2\theta+d)d+\theta+\frac{d}{2}-1)})\Pi_{\alpha,\theta,\nu_{0}}(|F|)^2, \quad r>0, F\in \mathcal{F}_{d}.
$$
\end{proof}

\subsection{Proof of Theorem \ref{d_{2}}}
\beg{proof}
If $d=\sharp\{i\in S, \nu_{0}(i)> 0\}=\infty,$
we follow the proof in \cite{S} to prove the invalidity of super Poincar\'e inequality for $\Pi_{\alpha,\theta,\nu_{0}}$ and $\EE.$
As defined before $p_i=\nu_{0}(i), i\geq 1,$ and $\lim_{i\rightarrow\infty}p_{i}=0.$
Let
$$F_{n}(\mu)=\bigg(\frac{1}{p_n(1+\theta p_n)}\bigg)^{\frac{1}{2}}\mu(n).$$
$\forall c> 0,$
\beg{align*}\label{keji}
&\int_{\{F_{n}^{2}\geq c\}}F_{n}^{2}d\Pi_{\alpha,\theta,\nu_{0}}
=\frac{1}{p_n(1+\theta p_n)}\int_{\{\mu|\mu(n)\geq\sqrt{cp_n(1+\theta p_n)}\}}\mu(n)^{2}\Pi_{\alpha,\theta,\nu_{0}}(d\mu)\\
&=\frac{1}{p_n(1+\theta p_n)}\int_{\sqrt{cp_n(1+\theta p_n)}}^{1}\frac{p_n(1-p_n)\Gamma(\theta+\frac{3}{2})}{\pi\Gamma(\theta+\frac{1}{2})}\frac{t^{\frac{1}{2}}(1-t)^{-\frac{3}{2}}}{(\frac{p_{n}^{2}}{t}+\frac{(1-p_{n})^{2}}{1-t})^{\theta+\frac{3}{2}}}dt\\
&\overrightarrow{n\rightarrow\infty}\frac{\Gamma(\theta+\frac{3}{2})}{\pi\Gamma(\theta+\frac{1}{2})}\int_{0}^{1}t^{\frac{1}{2}}(1-t)^{\theta}dt\\
&=\frac{\Gamma(\theta+\frac{3}{2})\Gamma(\frac{1}{2})\Gamma(\theta)}{\pi(\Gamma(\theta+\frac{1}{2}))^{2}}
,\end{align*}
which implies that $\{F_{n}^{2}\}_{n\geq 1}$ is not uniformly integrable.
So the F-Sobolev inequality doesn't hold, then neither will the super Poincar\'e inequality hold.
\end{proof}

\section*{Acknowledgement} The author would like to thank her doctoral supervisor Professor Feng-Yu Wang
for providing the idea of Theorem \ref{super Poincare 1},
and the work is supported in part by NNSFC (11771326).

\

\beg{thebibliography}{99}

\bibitem{Bakry} D. Bakry, I. Gentil, M. Ledoux, \emph{Analysis and geometry of markov diffusion operators}, Springer 2014.

\bibitem{MAC02} M.A.Carlton, \emph{A family of densities derived from the three-parameter Dirichlet process,} J. Appl. Prob (39)2002, 764-774.


 \bibitem{Davies} E. B. Davies, B. Simon, \emph{Ultracontractivity and the heat kernel
for Schr\"odinger operators and Dirichlet Laplacians,} J. Funct. Anal. 59(1984), 335--395.

\bibitem{QQ} C. L. Epstein, R. Mazzeo, \emph{Wright-Fisher diffusion in one dimension,}  SIAM J. Math. Anal.   42(2010),  568--608.






 \bibitem{FW14}  S. Feng,   F.-Y. Wang,   \emph{Harnack inequality and applications for infinite-dimensional GEM   processes,} Potential Anal.
  44(2016),  137--153

 \bibitem{FS18}  S. Feng,   W. Sun,   \emph{A dynamic model for the two-parameter Dirichlet process,} Potential Anal.
  1(2018),  18.

\bibitem{FSWX11}  S. Feng,   W. Sun, F-Y Wang, F. Xu,  \emph{Functional inequalities for the two-parameter extension of the infinitely-many-neutral-alleles diffusion,} Journal of functional analysis. 260(2011) 399-413.

\bibitem{Jac01} M. Jacobsen, \emph{Examples of multivariate diffusions: time-reversibility; a Cox-Ingersoll-Ross type process,} Department of Theoretical Statistics, Preprint 6, University of Copenhagen, 2001.

\bibitem{JCF} J. C. F. Kingman, \emph{Random discrete distributions,} J. Roy. statist. Soc. B. {\bf 37}(1975), 1-22.


 \bibitem{P2} J. E. Mosimann, \emph{On the compound multinomial distribution, the multivariate-distribution, and correlations among proportions,}
 Biometrika, 49(1962), 65--82.

 \bibitem{S} W. Stannat, \emph{On validity of the log-Sobolev
      inequality for symmetric Fleming-Viot operators,}  Ann.
      Probab.
      28(2000), 667--684.

\bibitem{W00a}
F.-Y. Wang, \emph{Functional inequalities for empty essential
spectrum, }  J. Funct. Anal. 170(2000), 219--245.

\bibitem{W00b} F.-Y. Wang, \emph{Functional inequalities, semigroup properties
and spectrum estimates,}   Infin. Dimens. Anal. Quant. Probab.
Relat. Topics 3(2000), 263--295.

\bibitem{Wbook} F.-Y. Wang,  \emph{Functional inequalities, markov semigroups and spectral theory,}   Science Press 2005.

\bibitem{WZ18} F.-Y. Wang, W.-W. Zhang,  \emph{Nash inequality for diffusion processes associated with Dirichlet distribution,} arXiv:1801.09209, 2018.

\end{thebibliography}
\end{document}